\documentclass{amsart}
\usepackage{amssymb,amsmath,amsfonts,amscd,euscript}

\newcommand{\nc}{\newcommand}

\nc{\slt}{\mathfrak{sl}_2}
\nc{\suth}{\widehat{\mathfrak{su}}(2)}
\nc{\gl}{\mathfrak{gl}}
\nc{\GL}{\mathfrak{GL}}
\nc{\g}{\mathfrak{g}}
\nc{\h}{\mathfrak{h}}
\nc{\la}{\mathfrak{a}}
\nc{\slth}{\widehat{\slt}}
\nc{\C}{\mathbb C }
\nc{\Z}{\mathbb Z }
\nc{\N}{\mathbb N }
\nc{\al}{\alpha }
\nc{\be}{\beta}
\nc{\ve}{\varepsilon}
\nc{\ch}{{\mathop {\rm ch}}}
\nc{\Tr}{{\mathop {\rm Tr}\,}}
\nc{\U}{{\mathop {\rm U}}}
\nc{\bra}{\langle}
\nc{\ket}{\rangle}
\nc{\x}{{\bf x}}
\nc{\pa}{\partial}
\nc{\ld}{\ldots}
\nc{\cd}{\cdots}
\nc{\sm}{\sqrt{2m}\,}
\nc{\hk}{\hookrightarrow}
\nc{\A}{\mathfrak A}
\nc{\qb}[2]{\genfrac{(}{)}{0pt}{}{#1}{#2}_q}
\nc{\n}{\mathfrak{n}}
\nc{\un}{\mathfrak{u}}
\nc{\T}{\otimes}
\nc{\bm}{{\bf m}}
\nc{\bs}{{\bf s}}
\nc{\bv}{{\bf v}}
\nc{\bt}{{\bf t}}
\nc{\bu}{{\bf u}}
\nc{\bin}[2]{{\genfrac{[}{]}{0pt}{0}{#1}{#2}}_q}
\nc{\fac}[1]{(#1)_q!}

\newtheorem{theo}{Theorem}[section]
\newtheorem{lem}{Lemma}[section]
\newtheorem{prop}{Proposition}[section]
\newtheorem{cor}{Corollary}[section]
\newtheorem{rem}{Remark}[section]

\begin{document}
\author{B.Feigin and E.Feigin}
\title
[Homological realization of restricted Kostka polynomials]
{Homological realization of restricted Kostka polynomials}

\address{Boris Feigin:
{\it Russia, Chernogolovka 142432, Landau Institute for The\-ore\-ti\-cal
Physics} and
{\it
Russia, Moscow, 119002, Independent University of Moscow, Bol'shoi Vlas'evskii, 11}}
\email{feigin@mccme.ru}
\address{Evgeny Feigin:
{\it Russia, Moscow, Moscow State University, Mech-Math Faculty,
Department of
Higher Algebra, Leninskie gori, 1} and {\it
Russia, Moscow, 119002, Independent University of Moscow, Bol'shoi Vlas'evskii, 11}}
\email{evgfeig@mccme.ru}

\begin{abstract}
In this paper we give two realizations of restricted Kostka polynomials
for $\slt$. Firstly we identify
 restricted Kostka polynomials with characters of the
zero homology of the current algebra with  coefficients
in the certain modules.
As a corollary we reobtain the alternating sum formula. Secondly we show that
restricted Kostka polynomials are  $q$-multiplicities of the
decomposition of the certain integrable $\slth$-modules to the irreducible
components. This allows to write a kind of fermionic formula for the Virasoro
unitary characters.
\end{abstract}
\maketitle

\section*{Introduction}
For $\bm\in\Z_{\ge 0}^k$ let $K_{l,\bm}(q)$ be a Kostka polynomial for
$\slt$ and $K^{(k)}_{l,\bm}(q)$ be a level
$k$ restricted Kostka polynomial
(we follow the notations in \cite{kos}; see also \cite{ss}).
Let $V_\bm$ be the fusion  product,
$$V_\bm=\underbrace{\pi_1*\cd *\pi_1 *}_{m_1}\cd *
\underbrace{\pi_k *\cd *\pi_k}_{m_k},$$
where $\pi_i$ is irreducible $(i+1)$-dimensional representation of $\slt$
(see \cite{fus, one, CP}).
Recall that $K_{l,\bm}(q)$ is the $q$-multiplicity of $\pi_l$ in
$V_\bm$ (see \cite{kos}). Therefore, $K_{l,\bm}(q)$ is equal to the
character of the $\slt$-invariants in $V_\bm\T\pi_l$.
We consider the induced module $\mathrm{Ind}_{\slt}^{\slt\T\C[t]} \pi_l$.
Note that $V_\bm\T \mathrm{Ind}_{\slt}^{\slt\T\C[t]} \pi_l$ is
a free $\slt\T t\C[t]$-module, and therefore, the relative homology
$$
H_0(\slt\T \C[t],\slt; V_\bm\T \mathrm{Ind}_{\slt}^{\slt\T\C[t]} \pi_l)
$$
is isomorphic to the $\slt$-invariants in the tensor product
$V_\bm\T \pi_l$. We obtain
\begin{equation}
\label{K}
K_{l,\bm}(q)=\ch_q H_0(\slt\T \C[t],\slt; V_\bm\T
\mathrm{Ind}_{\slt}^{\slt\T\C[t]} \pi_l).
\end{equation}

We now replace the induced module in $(\ref{K})$ by its quotient, which is
isomorphic to some irreducible integrable representation of $\slth$. Namely,
let $L_{l,k}$, $0\le l\le k$ be the set of
irreducible level $k$ integrable $\slth$-modules.
We show that
\begin{equation}
\label{0}
K^{(k)}_{l,\bm}(q)=\ch_q H_0(\slt\T\C[t],\slt; V_\bm\T L_{l,k}^*),
\end{equation}
and the higher homology vanish:
\begin{equation}
\label{p>0}
H_p(\slt\T\C[t], \slt; V_\bm\T L_{l,k}^*)=0,\ p>0.
\end{equation}
We note that $L_{l,k}^*=U(\slt\T\C[t])\cdot v_{l,k}^*$ ($v_{l,k}$ is a
highest
weight vector of $L_{l,k}$). Therefore, $L_{l,k}^*$ is a quotient of
$\mathrm{Ind}_{\slt}^{\slt\T\C[t]} \pi_l$.

Formulas $(\ref{0})$, $(\ref{p>0})$ and the BGG resolution allows us
to reobtain the alternating sum formula
\begin{equation*}
K^{(k)}_{l,\bm}(q)=\sum_{i\ge 0} q^{(k+2)i^2+(l+1)i}K_{2(k+2)i+l,\bm}(q)-
                \sum_{i> 0} q^{(k+2)i^2-(l+1)i}K_{2(k+2)i-l-2,\bm}(q).
\end{equation*}
The left-hand side is equal to
$\sum_{p\ge 0} (-1)^p \ch_q H_p(\slt\T\C[t],\slt; V_\bm\T L_{l,k})$,
while the right-hand side coincides with the Euler characteristics of
the certain complex, counting $H_*(\slt\T\C[t],\slt; V_\bm\T L_{l,k})$.

We now describe the realization  of restricted
Kostka polynomials as $q$-multi\-pli\-ci\-ties.
Let $L_{\bm,k}$ be an integrable  module induced from the fusion product
$V_\bm$:
\begin{equation}
\label{intind}
L_{\bm,k}=\left(\mathrm{Ind}_{\slt\T\C[t]}^{\slth} V_\bm\right)
/\bra K-k, e(z)^{k+1}\ket,
\end{equation}
where $K$ is the central element.
We show that this module coincides with the inductive limit of a fusion
products (see \cite{two}). Recall that
there exists an embedding $V_\bm\to V_{(m_1,\ld,m_{k-1}, m_k+2)}$.
We prove that $L_{\bm,k}=\lim_{s\to\infty} V_{(m_1,\ld,m_k+2s)}$.
Consider the decomposition of $L_{\bm,k}$ into the direct sum of
irreducible modules
\begin{equation}
\label{decompos}
L_{\bm,k}=\bigoplus_{l=0}^k  N_{l,\bm}\T L_{l,k}
\end{equation}
($N_{l,\bm}$ is spanned by highest weight vectors of the
weight $l$). It was proved in \cite{two} that the
following equality is true in the level $k$ Verlinde algebra
$V^{(k)}$
$$[1]^{m_1}[2]^{m_2}\cd [k]^{m_k}=[0]\dim N_{0,\bm} +\cd +
[k]\dim N_{k,\bm},$$
where $[l]$ corresponds to the $(l+1)$-dimensional representation of
$\slt$.
We show that the character of $N_{l,\bm}$ coincides with a reversed
Kostka polynomial
$\widetilde K^{(k)}_{l,\bm}(q)=q^{h(\bm)} K^{(k)}_{l,\bm}(q^{-1})$ for
the certain function $h:\Z_{\ge 0}^k\to \Z_{\ge 0}$.
This agrees with the fact that
$K^{(k)}_{l,\bm}(1)$ are the structure constants of $\mathrm{V}^{(k)}$ (see
\cite{HKKOTY,kos}).

We now apply the decomposition of $L_{\bm,k}$ to the coset construction
(see \cite{GKO}) to
obtain a finitization of the characters of the minimal unitary Virasoro
models. Our finitization is expressed in terms of Kostka polynomials
(see also \cite{ABF}, \cite{B}, \cite{Sch}).
We give some details below. (For
the connection of the branching functions and Kostka polynomials in a more
general settings see \cite{ss}. See also \cite{kedem, FKRW} for the
$\mathfrak{sl}_N$ case and $W_N$  instead of the Virasoro algebra).

The decomposition of $L_{\bm,k}$ can be applied to the study of the
coset constructions
$$(\slth)_{k_1}\T \cd \T (\slth)_{k_n}/(\slth)_{k_1+\cd +k_n}.$$
For example, we can get fermionic and bosonic formulas for the
corresponding conformal theories. In this paper we are dealing with
the simplest case $n=2$ and $k_1=1$.

Consider the decomposition of the tensor product
$L_{i,1}\T L_{j,k}=\bigoplus_{l=0}^{k+1} N_l\T L_{l,k+1}.$  Each $N_l$ is
a representation of the Virasoro algebra.
Namely, each $N_l$ is isomorphic to the certain minimal model
$M_{r,s}(k+2,k+3)$ (see \cite{KW}). We prove that
$L_{i,1}\T L_{j,k}=\lim_{N\to\infty} L_{(1^N (j+1)),k+1}$, where
$L_{(1^{m_1}\cd k^{m_k}),k}=L_{\bm,k}$.
This gives  a finitization of the characters of $N_l$ in terms of
Kostka polynomials.
We also show that the Rocha-Caridi formula
for the characters of the minimal models (see \cite{RC}) is a
corollary from the alternating sum formula.
There also exists the fermionic formula for restricted Kostka
polynomials (see \cite{ss,kos}). The limit of this formula
gives the fermionic type formula for the characters of the
minimal unitary  models (see \cite{B}). As mentioned in
\cite{JMT}, Kostka polynomials finitization for $j=0$ coincides with
the finitization in \cite{ABF}. We study this connection for general $j$ in
Appendix A.

We note that the fermionic formula for Kostka polynomials naturally
appears as $q$-multiplicity in the decomposition of $L_{\bm,k}$.
Namely, in \cite{two} the defining relations in principal subspaces
were described. This gives the fermionic formula for $L_{\bm,k}$.
We use a certain space of coinvariants to find the highest weight vectors
of the weight $l$ in the decomposition $(\ref{decompos})$.
This leads to the fermionic formula for the character of
$N_{l,\bm}$.

We finish the introduction with a discussion of possible generalizations.
Let $\widehat\g$ be an affine Kac-Moody algebra.
Fix $\widetilde \lambda=(\lambda_1,\ld,\lambda_k)$ to be a vector of the highest
weights
of the irreducible representations $\pi_{\lambda_i}$ of $\g$.
Let $V_{\widetilde \lambda}$ be the
corresponding fusion product $\pi_{\lambda_1}*\cd *\pi_{\lambda_k}$
(the adjoint graded space of the tensor product of the evaluation
representations, which is conjecturally independent of the evaluation
parameters, see \cite{fus, one, kedem, CL}).
Let $L_{\mu}$ be a level $k$ irreducible integrable highest weight representation
of $\widehat\g$ with highest weight $\mu$. Define
\begin{equation}
\label{gen1}
K^{(k)}_{\mu,{\widetilde \lambda}}(q)=
H_0(\g\T\C[t],\g;V_{\widetilde \lambda}\T L_{\mu}).
\end{equation}
(As above, $\g\T\C[t]$ are the generating operators in $L_\mu$).
We conjecture that
this definition of the restricted Kostka polynomials coincides with one in
\cite{kos}, as far as in \cite{SW}.

Another possibility is to use the induced module as in
$(\ref{intind})$. Namely, restricted Kostka polynomials can be
defined as the $q$-multiplicities of the
irreducible components in the decomposition of
\begin{equation}
\label{gen2}
\left(\mathrm{Ind}_{\g\T\C[t]}^{\widehat \g} V_{\widetilde \lambda}\right)
/\bra K-k, I\ket.
\end{equation}
Here $I\hk U(\widehat\g)$ is the ideal which vanishes in
any highest weight integrable level $k$ $\widehat\g$-module.

Our paper is organized as follows:

In Section $1$, we settle our notations and collect the main properties
of the fusion products.

In Section $2$, we give a homological realization of restricted Kostka
polynomials and prove the vanishing theorem for the higher homology.

In Section $3$, we derive the alternating sum formula from the BGG
resolution.

Section $4$ is devoted to the decomposition of $L_{\bm,k}$ into the
direct sum of irreducible modules.

In Section $5$, we obtain the finitization of the Virasoro unitary
characters in terms of Kostka polynomials and describe the coinvariants
approach to the fermionic formula.

In Appendix, the connection between the ABF and Kostka polynomials
finitizations of the unitary characters is studied.

{\bf Acknowledgment.} The first author was partially supported by the
RFBR grants 04-01-00303 and 05-01-01007, SS 2044.2003.2 and the 
INTAS grant 03-51-3350.
The second named author was partially supported by the RFBR grant
03-01-00167.

\section{Preliminaries}
\subsection{Fusion products}
In this section, we fix our notations and recall the main results in
\cite{fus,one, two}.

Let $e,h,f$ be the standard basis of $\slt$,
and let $\slth$ be the affine Kac-Moody algebra,
$$\slth=\slt\T\C[t,t^{-1}]\oplus \C K\oplus \C d,$$
where $K$ is a central element, and $[d,x_i]=-ix_i$.
We set $x_i=x\T t^i$ for $x\in\slt$, $i\in\Z$.

Let $\pi_l$ be $(l+1)$-dimensional irreducible
representation of $\slt$.
We fix some $n$-tuple $(z_1,\ld, z_n)$ of pairwise distinct complex numbers
and consider the tensor product
$$\pi_{a_1}(z_1)\T\cd \T\pi_{a_n}(z_n)$$
of the evaluation representations of $\slt\T\C[t]$. Then the fusion product
$$\pi_{a_1}*\cd *\pi_{a_n}$$
is an adjoint graded module with respect to the filtration
$$F_s=\mathrm{span}\bra x^{(1)}_{i_1}\cd x^{(t)}_{i_t} v_A,\
x^{(j)}\in\slt, i_1+\cd +i_t\le s,\ket$$
where $A=(a_1,\ld,a_n)$ and $v_A$ is the tensor product of the lowest weight
vectors of $\pi_{a_i}(z_i)$.

Suppose now that $a_i\le k$ for any $i$. Then we also use a following
notation. Let $\bm=(m_1,\ld, m_k)$ be some $k$-tuple of
non-negative integers. Set
$$V_\bm=\underbrace{\pi_1*\cd\* \pi_1}_{m_1} *\cd *
\underbrace{\pi_k *\cd *\pi_k}_{m_k}.$$
We recall that $V_\bm$ is cyclic with respect
to the algebra $\C[e_0, e_1, \ld]$. Denote the corresponding
cyclic vector by $v_\bm$ (note that this vector coincides with $v_A$).
We use a notation $u_\bm$ for a cyclic
vector with respect to the algebra $\C[f_0, f_1, \ld]$. Note that
$v_\bm$ ($u_\bm$) is the vector of the minimal (maximal) $h_0$-eigenvalue.
We recall that $V_\bm$ is bigraded, namely
\begin{equation}
\label{char}
V_{\bm}^{\al}=\{v\in V_{\bm}:\ h_0v=\al v\},\qquad
V_{\bm}^{\al,s}=V_{\bm}^{\al}\cap
\mathrm{span}\{e_{i_1}\cd e_{i_p}v_\bm:\ i_1+\cd +i_p=s\}.
\end{equation}
For $v\in V_{\bm}^{\al,s}$ set
$\deg_z v=\al,\ \deg_q v=s$, and define the character  of a homogeneous
subspace $V\hk V_\bm$ by
\begin{equation}
\label{ch}
\ch V=\sum_{\al,s} \dim (V\cap V_{\bm}^{\al,s}) z^\al q^s,\qquad
\ch_q V=\sum_{\al,s} \dim (V\cap V_{\bm}^{\al,s}) q^s.
\end{equation}
We now recall some exact sequences of fusion products.

Let $1\le a_1\le\cd\le a_n$. Then there exist
an exact sequences of $\slt\T\C[t]$-modules:
\begin{equation}
\label{our}
0\to \pi_{a_2-a_1}*\pi_{a_3}*\cd *\pi_{a_n}\to
\pi_{a_1}*\pi_{a_2}*\cd *\pi_{a_n}\to
\pi_{a_1-1}*\pi_{a_2+1}*\pi_{a_3}*\cd *\pi_{a_n}\to 0
\end{equation}
and
\begin{multline}
\label{tens}
0\to \pi_{a_1}*\cd *\pi_{a_{n-2}}\T \pi_{a_n-a_{n-1}}\to
\pi_{a_1}*\cd *\pi_{a_n}\to\\ \to
\pi_{a_1}*\cd *\pi_{a_{n-2}}*\pi_{a_{n-1}-1}*\pi_{a_n+1}\to 0.
\end{multline}

Now suppose that $a_i=a_{i+1}$. Then we also have an exact sequence of
$\slt\T\C[t]$-modules
\begin{multline}
\label{ind}
0\to \pi_{a_1}*\cd *\pi_{a_{i-1}} * \pi_{a_{i+2}} *\cd *\pi_{a_n}\to
\pi_{a_1}*\cd *\pi_{a_n}\to \\ \to
\pi_{a_1}*\cd *\pi_{a_{i-1}} * \pi_{a_i-1} *\pi_{a_{i+1}+1} *
\pi_{a_{i+2}} *\cd *\pi_{a_n}\to 0.
\end{multline}
We note that each of $(\ref{our}), (\ref{tens})$, and $(\ref{ind})$ contains
the piece
$$
\pi_{a_1}*\cd *\pi_{a_n}\to
\pi_{a_1}*\cd *\pi_{a_{i-1}} * \pi_{a_i-1} *\pi_{a_{i+1}+1} *
\pi_{a_{i+2}} *\cd *\pi_{a_n}
$$
in some special case. For the general case see \cite{three}.

We now recall some facts about a subspace
$$\C[e_1,e_2,\ld]\cdot v_A\hk \pi_{a_1} *\cd *\pi_{a_n}.$$
Note that $\C[e_1,e_2,\ld]\cdot v_A$ is
invariant with respect to  the subalgebra $\la_1\hk\slt\T\C[t]$
generated by $e_1$ and $f_0$.
Let $\la_2\hk\slt\T\C[t]$ be the subalgebra generated by $e_0$ and $f_1$.
Fix an isomorphism $\imath:\la_1\to\la_2$ sending $e_1$ to $e_0$ and
$f_0$ to $f_1$. Then we have an isomorphism of
$\la_2$-modules
\begin{equation}
\C[e_1,e_2,\ld]\cdot v_A\simeq \pi_{a_1}*\cd *\pi_{a_{n-1}},
\end{equation}
where the action of $\la_2$ on the left-hand side is a composition of
$\imath^{-1}$ and the natural action of $\la_1$.
In addition we have an exact sequence of $\la_1$-modules
\begin{equation}
\label{dem}
0\to \C[e_1,e_2,\ld]\cdot v_A\to
\pi_{a_1}*\cd *\pi_{a_n}\to \pi_{a_1}*\cd *\pi_{a_n-1}
\to 0.
\end{equation}

We describe an inductive limits of fusion products.
Using $(\ref{ind})$ we obtain a sequence of embeddings
\begin{equation}
\label{lim}
V_{(m_1,\ld,m_k)}\hk V_{(m_1,\ld,m_k+2)}\hk V_{(m_1,\ld,m_k+4)}\hk\cd.
\end{equation}
We denote the inductive limit of $(\ref{lim})$ by $L_{\bm,k}$.
This space can be endowed with the structure of
a level $k$ integrable $\slth$-module (the action
of the affine algebra is compatible with the natural action of the
annihilation operators $\slt\T\C[t]$). We consider the decomposition
\begin{equation}
\label{decomp}
L_{\bm,k}=L_{0,k}\T N_{0,\bm}\oplus \cd \oplus L_{k,k}\T N_{k,\bm},
\end{equation}
where $L_{i,k}$, $0\le i\le k$ are level $k$ irreducible highest weight
representations of $\slth$ with highest weight vectors $v_{i,k}$:
$h_0 v_{i,k}=iv_{i,k}$, $Kv_{i,k}=kv_{i,k}$, $dv_{i,k}=0$.
Then the dimensions of $N_{l,\bm}$ are given in terms of the level $k$
Verlinde algebra $\mathrm{V}^{(k)}$ for $\slt$.
Let $[0],[1],\ld,[k]$ be a basis of
$\mathrm{V}^{(k)}$ ($[l]$ corresponds to the $(l+1)$-dimensional
representation of $\slt$). Introduce the notation
\begin{equation}
[1]^{m_1}[2]^{m_2}\cd [k]^{m_k}=[0] c_{0,\bm} +\cd +
[k]c_{k,\bm}.
\end{equation}
The following theorem is proved in \cite{two}.
\begin{theo}
We have an equality
\begin{equation}
\label{Verl}
\dim N_{l,\bm}=c_{l,\bm}.
\end{equation}
\end{theo}

We finish this subsection with a remark on our characters notations.
For any homogeneous $V\hk L_{\bm,k}$ set
$$\ch V (z,q)=\sum_{\al,s} z^\al q^s \dim\{v\in V:\ h_0v=\al v, dv=sv\}  ,
\qquad  \ch_q V=\ch V (1,q).$$
We recall that the character $\ch_q V_\bm$ is given by $(\ref{ch})$.
We also need the "reversed" character, coming from the embedding
$\jmath: V_\bm\hk L_{\bm,k}$. Namely, set
$$\widetilde \ch_q V_\bm(q)=\ch_q (\jmath \cdot V_\bm).$$
Obviously, $\widetilde \ch_q V_\bm(q)=q^{h(\bm)}\ch_q V_\bm(q^{-1})$
for some $h: \Z_{\ge 0}^k\to\Z_{\ge 0}$ (see Lemma \ref{h(m)} for the
computation of $h(\bm)$).

\subsection{The Weyl group and Kostka polynomials}
We first settle our notations concerning $\slth$ (see \cite{Kac}).
Let $\h=\mathrm{span}\{h_0,K,d\}\hk\slth$ be the Cartan subalgebra,
$\n$ be the nilpotent subalgebra,
$\n=\slt\T t^{-1}\C[t^{-1}]\oplus\C f_0$, and
$\un=\slt\T t^{-1}\C[t^{-1}]$.
Let $s_0, s_1\in W$ be
simple reflections, where $W$ is the Weyl group of $\slth$.
Let $\rho\in\h^*$ be the element, defined by $\rho(\al_i^\vee)=1$, $i=0,1$,
where $\al_0^\vee$ and $\al_1^\vee$ are the simple coroots. We set
$w*\al=w(\al+\rho)-\rho$ for the shifted action of the Weyl group on
$\h^*$. Define $(i,k,m)\in\h^*$ by  $(i,k,m)h_0=i$,   $(i,k,m)K=k$, and
$(i,k,m)d=m$. Then $\rho=(1,2,0)$ and
\begin{equation}
\label{shift}
s_0*(i,k,m)=(-i+2k+2,k,m+k-i+1),\qquad s_1*(i,k,m)=(-i-2,k,m).
\end{equation}
The following lemma gives the shifted action of an arbitrary element of $W$
on $\h^*$.

\begin{lem}
\label{W}  $\phantom{a}$\\
$a).\quad s_0(s_1s_0)^n (i,k,m)=(-i-2+2(n+1)(k+2),k,
m+(n+1)^2(k+2)-(n+1)(i+1)).$\\
$b).\quad (s_0s_1)^n (i,k,m)=(i+2n(k+2),k,
m+n^2(k+2)+n(i+1)).$\\
$c).\quad s_1(s_0s_1)^n (i,k,m)=(-i-2-2n(k+2),k,
m+n^2(k+2)+n(i+1)).$\\
$d).\quad (s_1s_0)^n (i,k,m)=(i-2n(k+2),k,
m+n^2(k+2)-n(i+1)).$
\end{lem}

In our paper we use notations for Kostka polynomials as in
\cite{kos}. Let us recall the connection between the notations in
\cite{ss} and \cite{kos}.
For $\bm\in (\N\cup 0)^k$ we set
$$|\bm|=\sum_{i=1}^k im_i, \qquad 2\|\bm\|=-|\bm|+\sum_{1\le i,j\le k}
\min(i,j)m_im_j.$$
Now let $0\le l\le k$. Denote
$$\lambda=\left(\frac{|\bm|+l}{2}, \frac{|\bm|-l}{2}\right),\qquad
R(m)=(k^{m_k},\ld, 1^{m_1}).$$
We have
$$
K_{l,\bm}(q)=q^{\|\bm\|} K_{\lambda R(\bm)}(q^{-1}),\qquad
K^{(k)}_{l,\bm}(q)=q^{\|\bm\|} K^k_{\lambda R(\bm)}(q^{-1}),
$$
where the right-hand side stands for Kostka and level-restricted Kostka
polynomials in the notations of \cite{ss}.

\section{Homological realization of restricted Kostka polynomials}
In this section, we consider the fusion product $V_{\bm}$ as
$\slt\T\C[t^{-1}]$-module via the isomorphism
$$\slt\T\C[t^{-1}]\to \slt\T \C[t]\qquad x_i\mapsto x_{-i}.$$
Our goal is to show that
\begin{equation}
\label{main}
\ch_q H_p(\n,V_\bm\T L_{l,k})^0=\delta_{0,p} K^{(k)}_{l,\bm}(q),
\end{equation}
where $H_p(\n,V_\bm\T L_{l,k})^\al$ denotes an eigenspace of
the operator $h_0$ with an eigenvalue $\al$.
To prove this statement in the case $p=0$ we use the theorem in \cite{kos}:

\begin{theo}
\label{quot}
$\ch_q V_\bm/\bra h_0+l, e_0, e_{-1}^{k-l+1}\ket=K^{(k)}_{l,\bm}(q)$.
\end{theo}

\begin{lem}
\label{H_0}
$\ch_q H_0(\n,V_\bm\T L_{l,k})^0=K^{(k)}_{l,\bm}(q).$
\end{lem}
\begin{proof}
In view of Theorem $\ref{quot}$ it is enough to prove that
$$\ch_q H_0(\n,V_\bm\T L_{l,k})^0=
\ch_q V_\bm/\bra h_0+l, e_0, e_{-1}^{k-l+1}\ket.$$
We recall the first terms  of the BGG-resolution
\begin{equation}
\label{BGG}
0\gets L_{l,k}\gets M_{(l,k,0)}\xleftarrow{\partial}
M_{(-l-2,k,0)}\oplus M_{(-l+2k+2,k,k-l+1)}
\gets\cd,
\end{equation}
where $M_{\al}$ is the Verma module with a highest weight $\al$. We set
$v_{\al}$ to be a highest weight vector of $M_\al$.
Note that the differential $\partial$ is given by
\begin{equation}
\label{d}
\partial v_{(-l-2,k,0)}=f_0^{l+1} v_{(l,k,0)},\qquad
\partial v_{(-l+2k+2,k,k-l+1)}=e_{-1}^{k-l+1} v_{(l,k,0)}.
\end{equation}
We tensor $(\ref{BGG})$ by $V_\bm$ and obtain the free resolution of
the $\n$-module $L_{l,k}\T V_\bm$
\begin{equation*}
0\gets L_{l,k}\T V_\bm\gets M_{(l,k,0)}\T V_\bm \gets
(M_{(-l-2,k,0)} \oplus M_{(-l+2k+2,k,k-l+1)})\T V_\bm \gets\cd.
\end{equation*}
Therefore, $H_0(\n,V_\bm\T L_{l,k})^0$ is isomorphic to the homology
of the complex
\begin{multline*}
0\gets \left[\C\T _{U(\n)} \left(M_{(l,k,0)}\T V_\bm\right)\right]^0\gets\\
\gets \left[\C\T _{U(\n)}
\left( M_{(-l-2,k,0)}\T V_\bm\oplus M_{(-l+2k+2,k,k-l+1)}\T V_\bm\right)
\right]^0,
\end{multline*}
where $U(\n)$ is the universal enveloping algebra.
We note that $M_{(i,k,m)}$ is a free $U(\n)$-module with  one generator with
$z$-degree $i$. Hence, because of $(\ref{d})$
$$H_0(\n,V_\bm\T L_{l,k})^0\simeq V_\bm/
\bra f_0^{l+1}, e_{-1}^{k-l+1}, h_0+l\ket.$$
But
$$\left[ V_\bm/\bra f_0^{l+1}\ket \right]^{-l}\simeq
  \left[ V_\bm/\bra e_0\ket \right]^{-l}.$$
This finishes the proof of the lemma.
\end{proof}

Our next step is the proof of the statement $(\ref{main})$ in the case
$m_1+\cd +m_k=1$, i.e., when $V_\bm$ is a single representation $\pi_n$.
For this, we first recall the homology of $\un$.

\begin{lem}
\label{u}
$H_p(\un, L_{l,k})$ is isomorphic to
$\pi_{p(k+2)+l}$ for  even $p$ and to
$\pi_{p(k+2)+k-l}$ for odd $p$ as a representation of $\slt$.
\end{lem}

\begin{proof}
We recall (see \cite{hom},\cite{kum}) that
$H_p(\un, L_{l,k})$ is irreducible
$\slt$-module with the highest weight $(w_p*(l,k,0))(h_0)$,
where $w_p=\underbrace{s_0s_1s_0\cd}_p$. But Lemma $\ref{W}$ gives that
$(w_p*(l,k,0))(h_0)=p(k+2)+l+1$ for even $p$ and
$(w_p*(l,k,0))(h_0)=p(k+2)+k-l+1$ for odd $p$.
\end{proof}

\begin{prop}
\label{base}
$H_p(\n, \pi_n\T L_{l,k})^0$ is one-dimensional if
$$p \text{ is even and } n=p(k+2)+l \qquad { or } \qquad
  p \text{ is odd and } n=p(k+2)+k-l,$$
and vanishes otherwise.
\end{prop}

\begin{proof}
We note that $\un$ is an ideal in $\n$. Consider the
Hochschild-Serre spectral sequence (see \cite{CE}) with
$$E^2_{p,q}=H_p(\n/\un, H_q(\un,\pi_n\T L_{l,k})).$$
We first note that $\un$ acts trivially on $\pi_n$. Therefore,
$$H_q(\un,\pi_n\T L_{l,k})\simeq \pi_n\T H_q(\un, L_{l,k}).$$
In addition, $\n/\un$ is one-dimensional algebra $\C f_0$. We obtain that
\begin{multline}
\label{dim}
\dim H_p(\n, \pi_n\T L_{l,k})^0=
\dim \mathrm{coker}
\left(\left[\pi_n\T H_p(\un, L_{l,k})\right]^0   \xleftarrow{f_0}
\left[\pi_n\T H_p(\un, L_{l,k})\right]^2\right)+ \\
\dim \ker \left(\left[\pi_n\T H_{p-1}(\un, L_{l,k})\right]^0   \xleftarrow{f_0}
\left[\pi_n\T H_{p-1}(\un, L_{l,k})\right]^2\right).
\end{multline}
This gives that $\dim H_p(\n, \pi_n\T L_{l,k})^0$ vanishes unless
$n+1=\dim H_p(\un, L_{l,k})$, and in this case (because of Lemma $\ref{u}$)
the dimension is equal to $1$ for
$n=p(k+2)+l$ with even $p$ or $n=p(k+2)+k-l$ with odd $p$. Proposition
is proved.
\end{proof}

From the proof of the proposition we obtain a following corollary.
\begin{cor}
\label{-1}
$H_p(\n,\pi_n\T L_{l,k})^{-1}=0$ for any $p,n\ge 0$.
\end{cor}
\begin{proof}
We note that for any $n,n_1\ge 0$ the operator
$$f_0: \left[\pi_n\T\pi_{n_1}\right]^1\to
\left[\pi_n\T\pi_{n_1}\right]^{-1}$$
is an isomorphism. Now our corollary follows from the formula $(\ref{dim})$.
\end{proof}

\begin{cor}
\label{van}
For any  $\bm\in (\N\cup 0)^{k'}$ we have
$H_p(\n,V_\bm\T L_{l,k})^{-1}=0$, $p\ge 0$.
\end{cor}
\begin{proof}
Our corollary follows from Corollary $\ref{-1}$ and a fact that
$V_\bm$ has a filtration such that each quotient is
irreducible finite-dimensional $\slt$-module.
\end{proof}

\begin{rem}
We note that Corollary $\ref{van}$, $p=0$ follows from the exact
sequence $(\ref{BGG})$.
In fact, in the same way as in the proof of Lemma $\ref{H_0}$ we get
$$H_0(\n, V_\bm\T L_{l,k})^{-1}\simeq
V_\bm/\bra e_{-1}^{k-l+1}, f_0^{l+1}, h_0+l+1\ket.$$
But for any $v\in V_\bm\T L_{l,k}$ with $h_0v=-(l+1)v$ there exists
$v_1\in V_\bm\T L_{l,k}$ such that
$v=f_0^{l+1}v_1$. Therefore, $H_0(\n,V_\bm\T L_{l,k})^{-1}=0$.
\end{rem}

To prove that $H_p(\n,V_\bm\T L_{l,k})^0=0$ for $p>0$ we need one more
technical lemma. Let $\bm$ be some $k'$-tuple with $m_s\ne 0$ and
$m_{s+1}=\cd=m_{k'}=0$.
Introduce the notation for $k'$-tuple
$\bm^1=(m_1,\ld,m_s-1,0,\ld,0).$
We set
$$\widetilde V_{\bm^1}=U(\n)\cdot v_\bm\hk V_\bm$$
(recall that in this section fusion products are considered as modules
over the generating operators $x_i$, $i\le 0$).

\begin{lem}
\label{s}
$H_p(\n, \widetilde V_{\bm^1}\T L_{l,k})^\al\simeq
H_p(\n, V_{\bm^1}\T L_{k-l,k})^{k-s-\al}.$
\end{lem}

\begin{proof}
Consider the Lie algebra automorphism $\phi:\n\to\n$,
$\phi(e_{-1})=f_0$, $\phi(f_0)=e_{-1}$.
We note that $\phi$ induces an automorphism of the universal enveloping
algebra $U(\n)$. Denote this automorphism by the same letter.
Let $v_{l,k}$ be a
highest weight vector of $L_{l,k}$ and $v_\bm, u_\bm\in V_\bm$ be
lowest and highest (with respect to the operator $h_0$) weight vectors.
Then we have an isomorphism
$I:\widetilde V_{\bm^1}\T L_{l,k}\to V_{\bm^1}\T L_{k-l,k}$ defined by
$$x(v_\bm)\T y(v_{l,k})\mapsto
\phi(x)u_{\bm^1}\T \phi(y) v_{k-l,k},\quad x,y\in U(\n).$$
We note that $I$ is an isomorphism of $\n$-modules,
where the action of $\n$ on $V_{\bm^1}\T L_{k-l,k}$ is a composition
of $\phi$ and a standard action.

We need to show that $I$ identifies $(\widetilde V_{\bm^1}\T L_{l,k})^\al$ and
$(V_{\bm^1}\T L_{k-l,k})^{k-s-\al}.$ Note that
$[h_0,\phi(f_0)]=2\phi(f_0)$ and $[h_0,\phi(e_{-1})]=-2\phi(e_{-1})$.
Therefore, for $x\in U(\n)$ with $[h_0,x]=\be x$ ($\be\in\C$) we have
$[h_0,\phi(x)]=-\be \phi(x)$. In addition, in view of
$h_0 v_\bm=(-\sum_{i=1}^{k'} im_i)v_\bm$ and
$h_0 u_\bm=(\sum_{i=1}^{k'} im_i)u_\bm$, we obtain
\begin{multline}
h_0(v_\bm\T v_{l,k})=(l-\sum_{i=1}^{k'} im_i)v_\bm\T v_{l,k},\\
\shoveleft{h_0(I(v_\bm\T v_{l,k}))=h_0(u_{\bm^1}\T v_{k-l,k})=}\\
(\sum_{i=1}^{k'} im^1_i+k-l)v_\bm\T v_{l,k}=
(-l+\sum_{i=1}^{k'} im_i+k-s)v_\bm\T v_{l,k}.
\end{multline}
This finishes the proof of the lemma.
\end{proof}

\begin{cor}
\label{k+2}
Let $\bm\in \Z_{\ge 0}^{k+1}$, $m_{k+1}\ne 0$. Then
$H_p(\n, V_\bm\T L_{l,k})^0=0$ for any $p\ge 0$.
\end{cor}
\begin{proof}
We set $\bm^1=(m_1,\ld,m_{k+1}-1)$,  $\bm^2=(m_1,\ld,m_k+1,m_{k+1}-1)$.
Recall the exact sequence $(\ref{dem})$ of $\n$-modules
$$0\to \widetilde V_{\bm^1}\to V_\bm\to V_{\bm^2}\to 0.$$
We note that the map $V_\bm\to V_{\bm^2}$ is defined by
$u_\bm\mapsto u_{\bm^2}$ (because $V_\bm$ is cyclic $\n$-module with the
cyclic vector $u_\bm$). In addition,
$\deg_z u_{\bm^2}=\deg_z u_\bm-1$. Therefore, for any $\al$ we obtain an
exact sequence
$$0\to \widetilde V_{\bm^1}^\al\to V_\bm^\al\to V_{\bm^2}^{\al-1}\to 0.$$
This gives an exact sequence of homology
\begin{multline}
\label{al-1}
0\gets H_0(\n, V_{\bm^2}\T L_{l,k})^{\al-1}\gets
H_0(\n, V_\bm\T L_{l,k})^\al\gets\\
H_0(\n, \widetilde V_{\bm^1}\T L_{l,k})^\al\gets
H_1(\n, V_{\bm^2}\T L_{l,k})^{\al-1}\gets\cd
\end{multline}
Now let $\al=0$. Then because of Lemma $\ref{s}$ and Corollary $\ref{van}$
we obtain
$$H_p(\n, V_{\bm^2}\T L_{l,k})^{-1}=0, \quad
H_p(\n, \widetilde V_{\bm^1}\T L_{l,k})^0\simeq
H_p(\n, V_{\bm^1}\T L_{k-l,k})^{k-k-1}=0.$$
In view of the exact sequence $(\ref{al-1})$ our corollary is proved.
\end{proof}

We now prove the main theorem of this section.
\begin{theo}
\label{mth}
Let $\bm\in\Z_{\ge 0}^k$. Then
$$\ch_q H_p(\n,V_\bm\T L_{l,k})^0=\delta_{0,p} K^{(k)}_{l,\bm}(q).$$
\end{theo}
\begin{proof}
Because of Lemma $\ref{H_0}$ we only need to prove that
$H_p(\n,V_\bm\T L_{l,k})^0$ vanishes for $p>0$. We use the induction on $\bm$.
We order a $k$-tuples by the rule
$$\bm> {\bf n}\quad \text{ if } \sum_{i=1}^k m_i>\sum_{i=1}^k n_i
\text{ or }
\sum_{i=1}^k m_i=\sum_{i=1}^k n_i \text{ and }
\prod_{i=1}^k i^{m_i}>\sum_{i=1}^k i^{n_i}.$$
For $\sum_{i=1}^k m_i=1$ our theorem follows from Proposition $\ref{base}$.
Now let $V_\bm=\pi_{a_1}*\cd * \pi_{a_n}$ and $1\le a_1 \le \cd\le a_n$.
We recall the exact sequence $(\ref{our})$
\begin{equation}
\label{ind1}
0\to \pi_{a_2-a_1}*\pi_{a_3}*\cd *\pi_{a_n}\to
\pi_{a_1}*\pi_{a_2}*\cd *\pi_{a_n}\to
\pi_{a_1-1}*\pi_{a_2+1}*\pi_{a_3}*\cd *\pi_{a_n}\to 0.
\end{equation}
We denote the first fusion product in $(\ref{ind1})$ by $V_{\bm^{(1)}}$
and the third one by $V_{\bm^{(2)}}$. From $(\ref{ind1})$ we obtain
a long exact sequence
\begin{multline}
\label{long}
0\gets H_0(\n,V_{\bm^{(2)}}\T L_{l,k})\gets
H_0(\n,V_\bm\T L_{l,k})\gets H_0(\n,V_{\bm^{(1)}}\T L_{l,k})\gets\\
H_1(\n,V_{\bm^{(2)}}\T L_{l,k})\gets
H_1(\n,V_\bm\T L_{l,k})\gets H_1(\n,V_{\bm^{(1)}}\T L_{l,k})\gets\cd
\end{multline}
We note that $\bm^{(1)}< \bm$, and therefore,
$H_p(\n,V_{\bm^{(1)}}\T L_{l,k})^0$ vanishes for $p>0$ by induction assumption.
In addition, $\bm^{(2)}< \bm$, and
$V_{\bm^{(2)}}$ is either the fusion product of the representations of
dimension at most $k+1$ or one of the fused representations is of
the dimension $k+2$. In the latter case  Corollary $\ref{k+2}$ gives the vanishing
of the higher homology. Hence, because of the exact sequence
$(\ref{long})$, the theorem is proved.
\end{proof}

\begin{cor}
We have
$$\ch_q H_p(\slt\T\C[t^{-1}], \slt; V_\bm\T L_{l,k})=
\delta_{0,p} K^{(k)}_{l,\bm}(q).$$
\end{cor}
\begin{proof}
We note that as $\slt$-module $V_\bm\T L_{l,k}$ decomposes into the direct
sum of finite-dimensional representations. Therefore,
$$\ch_q H_p(\slt\T\C[t^{-1}], \slt; V_\bm\T L_{l,k})=
 \ch_q (H_p(\un, V_\bm\T L_{l,k})^{\slt})$$
(the right-hand side is a subspace of $\slt$-invariants).
In view of the Hochschild-Serre spectral sequence we also obtain
$$H_p(\n,V_\bm\T L_{l,k})^0=H_0(\C f_0,H_p(\un,V_\bm\T L_{l,k}))^0=
H_p(\un, V_\bm\T L_{l,k})^{\slt}.$$
Corollary is proved.
\end{proof}

\begin{rem}
We note that our theorem concerns only the case of a homology with the
coefficients in $\pi_{a_1}*\cd *\pi_{a_n}\T L_{l,k}$ with
$a_i\le k$. For the general $a_i$ the corresponding homology are not
concentrated in one dimension.
\end{rem}

\section{The BGG resolution and alternating sum formula for
Kostka polynomials}
In this section, we give a homological interpretation of the alternating
sum formula (see \cite{ss,kos})
\begin{equation}
\label{asf}
K^{(k)}_{l,\bm}(q)=\sum_{i\ge 0} q^{(k+2)i^2+(l+1)i}K_{2(k+2)i+l,\bm}(q)-
                \sum_{i> 0} q^{(k+2)i^2-(l+1)i}K_{2(k+2)i-l-2,\bm}(q).
\end{equation}
We note that in view of Theorem $\ref{mth}$ the left-hand side coincides
with the Euler characteristics
$\sum_{p\ge 0} (-1)^p \ch_q H_p(\n,V_\bm\T L_{l,k})^0$. The idea is that
the right-hand side is also the Euler characteristics of the complex,
counting the homology $H_p(\n,V_\bm\T L_{l,k})^0$.

Consider the BGG-resolution of $L_{l,k}$  (see \cite{BGG, kum})
\begin{equation}
\label{fBGG}
0\gets L_{l,k}\gets F_0\gets F_1\gets\cd,\quad
F_p=\bigoplus_{l(w)=p} M(w*(l,k,0)),
\end{equation}
where $l(w)$ is a length of the element of the Weyl group of $\slth$.
Tensoring $(\ref{fBGG})$ with $V_\bm$ we obtain the resolution for
$V_\bm\T L_{l,k}$.
Therefore, the following complex counts $H_*(\n,V_\bm\T L_{l,k})^0$
\begin{equation}
\label{Fp}
0\gets (\C\T_{U(\n)} (F_0\T V_\bm))^0\gets
(\C\T_{U(\n)} (F_1\T V_\bm))^0\gets\cd
\end{equation}
Recall that $F_p$ is a free $U(\n)$-module. Therefore, we can rewrite
$(\ref{Fp})$ as
\begin{equation*}
0\gets ((\C\T_{U(\n)} F_0)\T V_\bm)^0\gets
((\C\T_{U(\n)} F_1)\T V_\bm)^0\gets\cd
\end{equation*}

\begin{lem}
\label{chFp}
$\ch_q ((\C\T_{U(\n)} F_p)\T V_\bm)^0=
\sum_{l(w)=p} q^{(w*(l,k,0))d} \ch_q V_\bm^{-(w*(l,k,0))h_0}$.
\end{lem}
\begin{proof}
$F_p$ is a free $U(\n)$-modules with a generators labeled
by $w$ such that $l(w)=p$. In addition, the $z$-degree of the generator equals to
$(w*(l,k,0))h_0$. This proves our lemma.
\end{proof}

\begin{lem}
$K_{l,\bm}(q)=\ch_q V_\bm/\bra e_0^{l+1}, h_0+l \ket$, i.e.,
$K_{l,\bm}(q)$ is a $q$-multiplicity of $\pi_l$ in decomposition of $V_\bm$
as $\slt$-module to the irreducible components.
\end{lem}
\begin{proof}
Follows from Theorem $\ref{quot}$ and a fact
$\lim_{k\to\infty} K^{(k)}_{l,\bm}(q)=K_{l,\bm}(q)$.
\end{proof}

\begin{cor}
\label{fk}
$K_{l,\bm}(q)=\ch_q V_\bm^l-\ch_q V_\bm^{l+2}.$
\end{cor}

\begin{prop}
\begin{multline}
\label{Eul}
\sum_{p\ge 0} (-1)^p \ch_q (\C\T_{U(\n)} (F_p\T V_\bm))^0=\\
\sum_{p\ge 0} q^{(k+2)p^2+(l+1)p}K_{2(k+2)p+l,\bm}-
                \sum_{p> 0} q^{(k+2)p^2-(l+1)p}K_{2(k+2)p-l-2,\bm}.
\end{multline}
\end{prop}
\begin{proof}
In view of Lemma $\ref{chFp}$ and Lemma $\ref{W}$ we obtain that the left-hand
side of $(\ref{Eul})$ is equal to
\begin{multline*}
\sum_{p\ge 0}
\Bigl( q^{(p+1)^2(k+2)-(p+1)(l+1)} \ch_q (V_\bm)^{2p(k+2)-l}+ \\
         q^{p^2(k+2)+p(l+1)} \ch_q (V_\bm)^{-2p(k+2)-l}\Bigr)-\\
\shoveleft
{\sum_{p\ge 0}
\Bigl( q^{(p+1)^2(k+2)-(p+1)(l+1)} \ch_q (V_\bm)^{-(2p+1)(k+2)+l-k}+}\\
         q^{p^2(k+2)+p(l+1)} \ch_q (V_\bm)^{(2p+1)(k+2)-k+l}\Bigr)=\\
\shoveleft{\sum_{p\ge 0}
q^{p^2(k+2)+p(l+1)}
\left(\ch_q (V_\bm)^{2p(k+2)+l} -\ch_q (V_\bm)^{2p(k+2)+l+2}\right)+}\\
\sum_{p\ge 1}
q^{p^2(k+2)-p(l+1)}
\left(\ch_q (V_\bm)^{2p(k+2)-l}-\ch_q (V_\bm)^{2p(k+2)-l-2}\right)=\\
\sum_{p\ge 0}
q^{p^2(k+2)+p(l+1)} K_{2p(k+2)+l,\bm}(q)-
\sum_{p\ge 1}  q^{p^2(k+2)-p(l+1)} K_{2p(k+2)-l-2,\bm}(q),
\end{multline*}
where Corollary $\ref{fk}$ is used.
\end{proof}

As a corollary we obtain the alternating sum formula $(\ref{asf})$.

\section{The decomposition of $L_{\bm,k}$}
We introduce the notations $\widetilde K_{l,\bm}(q)$ and
$\widetilde K^{(k)}_{l,\bm}(q)$ for the "reversed"
Kostka polynomials:
\begin{equation}
\label{rev}
\widetilde K^{(k)}_{l,\bm}(q)=q^{h(\bm)} K^{(k)}_{l,\bm}(q^{-1}),\quad
\widetilde K_{l,\bm}(q)=q^{h(\bm)} K_{l,\bm}(q^{-1}),
\end{equation}
where $h(\bm)=\max\{ \deg_q v,\ v\in V_\bm\}$
(see Lemma $\ref{h(m)}$ for the computation of $h(\bm)$).
Therefore, "reversed' polynomials are  Kostka polynomials in
notations in \cite{ss} up to a power of $q$.

We recall the decomposition $(\ref{decomp})$
$$L_{\bm,k}=L_{0,k}\T N_{0,\bm}\oplus \cd \oplus L_{k,k}\T N_{k,\bm}.$$
\begin{lem}
\label{rK}
$\ch_q N_{l,\bm}=\widetilde K^{(k)}_{l,\bm}(q)$.
\end{lem}
\begin{proof}
We first note that
$K^{(k)}_{l,\bm}(1)=\dim N_{l,\bm}$, because both sides are
structure constants of the Verlinde algebra (see $(\ref{Verl})$).

We show that $\ch_q N_{l,\bm}=\widetilde
\ch_q V_\bm/\bra f_0, h_0-l, f_1^{k-l+1}\ket$. Recall that
$L_{\bm,k}=\lim_{s\to\infty} V_{\bm(s)}$, where
$\bm(s)=(m_1,\ld,m_{k-1},m_k+2s)$.
We denote
$$N^s_{l,\bm}=\{v\in L_{\bm,k}:\ e_0v=f_1v=0, h_0v=lv,
v\in V_{\bm(s)}\}.$$
Note that $N^s_{l,\bm}$ is a subspace of the space of highest weight
vectors of the weight $l$. Therefore, any $v\in N^s_{l,\bm}$ is not the
element of  $\bra f_0, h_0-l, f_1^{k-l+1}\ket V_\bm$ (because
for a highest weight vector $v_{l,k}\in L_{l,k}$ we have
$v_{l,k}\notin \bra f_0, h_0-l, f_1^{k-l+1}\ket L_{l,k}$).
We thus obtain that
$$\ch_q N^s_{l,\bm}\le \ch_q V_{\bm(s)}/\bra f_0, h_0-l, f_1^{k-l+1}\ket$$
(the difference of the right-hand side and left-hand side is a polynomial
with  nonnegative coefficients).
In addition, there exists $s_0$ such that $N^{s_0}_{l,\bm}=N_{l,\bm}$.
Hence, $\ch_q N_{l,\bm}\le \widetilde K^{(k)}_{l,\bm(s_0)}(q)$.
Let $c_{l,\bm}$ denote the structure constants of the level $k$
Verlinde algebra:
$$[1]^{m_1}\cd [k]^{m_k}=\sum_{l=0}^k  c_{l,\bm}[l].$$
In view of $[k]^2=[0]$ we get
$\ch_q N_{l,\bm}(1)=c_{l,\bm}=\widetilde K^{(k)}_{l,\bm(s_0)}(1)$
(for the second equality see \cite{HKKOTY}, \cite{kos}).  Therefore,
there exists $s_0$ such that
$$\ch_q N_{l,\bm}= \widetilde K^{(k)}_{l,\bm(s_0)}(q).$$
To complete the proof we need to show that
\begin{equation}
\label{indep}
\widetilde K^{(k)}_{l,\bm}(q)=\widetilde K^{(k)}_{l,\bm(1)}(q).
\end{equation}

We recall an exact sequence of $\slt\T\C[t]$-modules
$$0\to V_\bm\to V_{\bm(1)}\to
V_{(m_1,\ld,m_{k-2}, m_{k-1}+1,m_k,1)}\to 0.$$
Because of Corollary $\ref{k+2}$ the corresponding long exact sequence of
$\n$-homology is of the form
$$0\gets 0\gets H_0(\n, V_{\bm(1)}\T L_{l,k})^0\gets
H_0(\n, V_{\bm}\T L_{l,k})^0 \gets 0\gets\cd.$$
Because of Lemma $\ref{H_0}$ the equation $(\ref{indep})$ is shown.
\end{proof}

We finish this section with the identification of $L_{\bm,k}$ with the
induced module from the fusion product $V_\bm$. We first need one lemma.

\begin{lem}
\label{=0}
Let $1\le a_1\le\cd\le a_n$ and $a_n\ge k+1$. Then
$\pi_{a_1}*\cd \pi_{a_n}/\bra e_0, h_0+l, e_1^{k-l+1}\ket=0$ for
any $0\le l\le k$.
\end{lem}
\begin{proof}
We prove our lemma by induction on a pair $(n,\prod_{i=1}^n a_i)$.
We set $(n_1,s_1)>(n_2,s_2)$ if $n_1>n_2$ or $n_1=n_2$ and $s_1>s_2$.
Let $n=1$. Then for any $0\le l\le k$ we have
$\pi_a/\bra h_0+l, e_0\ket=0$ if $a\ge k+1$. We now consider an exact
sequence of $\slt\T \C[t]$-modules
$$0\to \pi_{a_2-a_1}*\pi_{a_3}*\cd *\pi_{a_n}\to
\pi_{a_1}*\cd *\pi_{a_n}\to
\pi_{a_1-1}*\pi_{a_2+1}*\pi_{a_3}*\cd *\pi_{a_n}\to 0.
$$
By induction assumption our lemma is true for the submodule and for the
quotient module. Therefore, it also holds for $\pi_{a_1}*\cd *\pi_{a_n}$.
\end{proof}

\begin{prop}
Let $\bm\in\Z_{\ge 0}^k$. Then
\begin{equation}
\label{induced}
L_{\bm,k}=\left(\mathrm{Ind}_{\slt\T\C[t]}^{\slth} V_\bm\right)
/\bra e(z)^{k+1}, K-k\ket,
\end{equation}
where the right-hand side is a quotient of the induced module
(with fixed $K=k$) by the action of the coefficients of the series
$e(z)^{k+1}=(\sum_i e_i z^i)^{k+1}$. In addition
\begin{equation}
\label{vinduced}
\left(\mathrm{Ind}_{\slt\T\C[t]}^{\slth} \pi_{a_1}*\cd *\pi_{a_n}\right)
/\bra e(z)^{k+1}, K-k \ket=0,
\end{equation}
if $a_1\le\cd\le a_n$ and $a_n\ge k+1$.
\end{prop}
\begin{proof}
Note that the right-hand side of $(\ref{induced})$ is a level $k$
integrable $\slth$-module (because $e(z)^k$ acts by $0$). Therefore, it
can be decomposed into the direct sum of irreducible modules $L_{l,k}$.
We show $(\ref{induced})$ by checking that the $q$-multiplicity of
$L_{l,k}$ is equal to
$\widetilde \ch_q V_\bm/\bra e_0, h_0+l, e_1^{k-l+1}\ket$.

Consider $\slth$ homomorphisms  from $L_{l,k}$ to the right-hand side of
$(\ref{induced})$.  They coincide with the homomorphisms
$L_{l,k}\to \mathrm{Ind}_{\slt\T\C[t]}^{\slth} V_\bm/\bra K-k\ket$,
which are labeled by the elements of the quotient
$V_\bm/\bra e_0, h_0+l, e_1^{k-l+1}\ket$ (because of the highest weight
condition for $L_{l,k}$). The first part of our proposition is verified.

We note that $(\ref{vinduced})$ can be checked in the same manner, taking
into account Lemma $\ref{=0}$.
\end{proof}

\section{Virasoro unitary models}
\subsection{The alternating sum formula.}
We first recall the coset construction (see \cite{GKO}).
Consider the decomposition of the tensor product $L_{i,1}\T L_{j,k}$
into the sum of irreducible $\slth$-modules
$$L_{i,1}\T L_{j,k}=\bigoplus_{l=0}^{k+1} N_l\T L_{l,k+1},$$
where $N_l$ is spanned by a highest weight vectors of
the weight $l$.
Let
$L^{(1)}_i, L^{(2)}_i$ and $L^{diag}_i$ be the Sugawara operators, acting
on $L_{i,1}$, $L_{j,k}$ and $L_{i,1}\T L_{j,k}$. Then
operators
$$L_i=L^{(1)}\T \mathrm{Id}+\mathrm{Id}\T L^{(2)}_i-L^{diag}_i$$ form the
Virasoro algebra, which acts on the tensor product $L_{i,1}\T L_{j,k}$
with the central charge $\frac{3}{3+2}+\frac{3k}{k+2}-\frac{3(k+1)}{k+3}=
\frac{k^2+5k}{(k+2)(k+3)}$.
The important property is that $L_i$ commute with the diagonal action
of $\slth$. Therefore, each $N_l$ is a representation of the Virasoro
algebra. Using the alternating sum formula we
derive a formula for the character of $N_l$. This formula
coincides with the Rocha-Caridi formula for the character of the minimal
model $M_{j+1,l+1}(k+2,k+3)$ (see \cite{KW}).

We first recall the Rocha-Caridi formula for the character of
$M_{r,s}(p,p')$ (see \cite{RC}).
Here $p,p'$ are relatively prime numbers and $1\le r\le p-1$,
$1\le s\le p'-1$.
Let $t=\frac{p'}{p}$. Then the central element $c$ of the Virasoro
algebra acts on $M_{r,s}(p,p')$ as a scalar $13-6(t+\frac{1}{t})$.
Let $\triangle_{r,s}=\frac{(rt-s)^2-(t-1)^2}{4t}$ and
$$\chi_{r,s}=\ch M_{r,s}(p,p')=\sum_{d\in \Z_{\ge 0}+\triangle_{r,s}} q^d
\dim \{v:\ L_0v=dv\}.$$
We set
$$ (q)_n=\prod_{\al=1}^n (1-q^\al),\quad
(q)_\infty=\prod_{\al=1}^\infty (1-q^\al),\quad
\bin{m}{n}=\frac{(q)_m}{(q)_n (q)_{m-n}}.
$$
Then
$$\chi_{r,s}=\frac{q^{\triangle_{r,s}}}{(q)_\infty}
\left(\sum_{n\in\Z} q^{pp'n^2+(p'r-ps)n}-
\sum_{n\in\Z} q^{pp'n^2+(p'r+ps)n+rs}\right).$$
We note that in the case $(p,p')=(k+2,k+3)$ the central charge is
equal to $\frac{k(k+5)}{(k+2)(k+3)}$.

We now recall the embedding of $V_\bm$ and $L_{\bm,k}$ into the tensor product
of the level one irreducible modules.
In what follows we use the notation $L_{(1^{m_1}\cd k^{m_k}),k}$ for
$L_{\bm,k}$.
Let $v(p)\in L_{0,1}\oplus L_{1,1}$
be the set of extremal vectors, $h_0v(p)=-pv(p)$. Then
we have the isomorphisms
\begin{gather}
\label{embV}
V_\bm\simeq
U(\slt\T \C[t])\cdot (v(m_1+\cd +m_k)\T v(m_2+\cd +m_k)\T\cd\T v(m_k)),\\
\label{embL}
L_{\bm,k}\simeq
U(\slth)\cdot (v(m_1+\cd +m_k)\T v(m_2+\cd +m_k)\T\cd\T v(m_k)).
\end{gather}
This gives an embedding of $V_\bm$ and $L_{\bm,k}$ into the tensor product
$L_{i_1,1}\T\cd\T L_{i_k,1}$,
where $i_\al=0,1$.  For example, in view of
$L_{j,k}=L_{(j^1),k}$ we obtain
\begin{equation*}
L_{j,k}\simeq U(\slth)\cdot (v(1)^{\T j}\T v(0)^{\T (k-j)})\hk
L_{1,1}^{\T j}\T L_{0,1}^{\T (k-j)}.
\end{equation*}

\begin{lem}
We have an isomorphism of $\slth$-modules
\begin{equation}
\label{limit}
L_{i,1}\T L_{j,k}\simeq\lim_{n\to\infty}
U(\slth)\cdot (v(2n+i)\T v(1)^{\T j}\T v(0)^{\T (k-j)}).
\end{equation}
\end{lem}
\begin{proof}
We first note that $v(2n+i)\T [v(1)^{\T j}\T v(0)^{\T (k-j)}]$ is the tensor
product of extremal vectors of $L_{i,1}$ and $L_{j,k}$. In addition,
\begin{multline*}
U(\slth)\cdot (v(2n+i)\T v(1)^{\T j}\T v(0)^{\T (k-j)})=\\
U(\slth)\cdot (v(2n+i+2s)\T v(2s+1)^{\T j}\T v(2s)^{\T (k-j)})
\end{multline*}
for any integer $s$. Therefore, to prove our lemma it suffices to show that
$$v(2n+i)\T v(1)^{\T j}\T v(0)^{\T (k-j)}\hk
U(\slth)\cdot (v(2n+2+i)\T v(1)^{\T j}\T v(0)^{\T (k-j)})$$
(in this case all the products of the extremal vectors of $L_{i,1}$ and
$L_{j,k}$ are the elements of $(\ref{limit})$).

We recall that $e_{N-1}v(N)=v(N-2)$ and $e_{\ge N}v(N)=0$. Hence,
$$e_{i+2n-1} (v(2n+i)\T v(1)^{\T j}\T v(0)^{\T (k-j)})=
(v(2n-2+i)\T v(1)^{\T j}\T v(0)^{\T (k-j)}).$$
Lemma is proved.
\end{proof}

\begin{cor}
$L_{i,1}\T L_{j,k}\simeq\lim_{n\to\infty} L_{(1^{2n+i-1} (j+1)),k+1}.$
\end{cor}
\begin{proof}  Because of the formula $(\ref{embL})$
$$U(\slth)\cdot (v(2n+i)\T v(1)^{\T j}\T v(0)^{\T (k-j)})\simeq
L_{(1^{2n+i-1} (j+1)),k+1}.$$
\end{proof}

Consider the decomposition
$$L_{(1^N (j+1)),k+1}=\bigoplus_{l=0}^{k+1} N_{l,(1^N (j+1))}\T L_{l,k+1}.$$

\begin{cor}
We have
$$\ch_q N_l=\lim_{n\to\infty} \ch_q N_{l,(1^{2n+i-1} (j+1))}.$$
\end{cor}

We now compute the limit from the above corollary.
Because of Lemma $\ref{rK}$
$$\ch_q N_{l,(1^N (j+1))}=\widetilde K^{(k+1)}_{l,(1^N (j+1))}(q).$$

Recall that $L_0 v_{i,k}=\frac{i(i+2)}{4(k+2)} v_{i,k}$ for the highest
weight vector $v_{i,k}\in L_{i,k}$.
In view of the alternating sum formula $(\ref{asf})$ and the formula
$(\ref{rev})$ we obtain
\begin{multline}
\label{chN}
\ch_q N_{l,(1^{2n+i-1} (j+1))}=q^{\frac{i(i+2)}{12}+\frac{j(j+2)}{4(k+2)}-
\frac{l(l+2)}{4(k+3)}} \widetilde K^{(k+1)}_{l,(1^{2n+i-1} (j+1))}(q)=\\
=q^{\frac{i(i+2)}{12}+\frac{j(j+2)}{4(k+2)}-\frac{l(l+2)}{4(k+3)}}
\biggl(\sum_{p\ge 0} q^{-(k+3)p^2-(l+1)p}
\widetilde K_{2(k+3)p+l,(1^{2n+i-1} (j+1))}(q)-\\
\sum_{p> 0} q^{-(k+3)p^2+(l+1)p}
\widetilde K_{2(k+3)p-l-2,(1^{2n+i-1} (j+1))}(q)\biggr)
\end{multline}
We want to compute the limit of the above expression while $n\to\infty$.

\begin{lem}
Let $a(q)=\sum_{i\ge 0} a_iq^i$. We write $a(q)=\mathrm{O}(q^N)$ if $a_i=0$
for $i<N$. Then
\begin{multline*}
\widetilde K_{2s+i+j,(1^{2n+i+1} (j+1))}(q)-
(\ch_q L_{i,1}^{2s+i}-\ch_q L_{i,1}^{2s+i+2j+2})=\mathrm{O}(q^{n+s(s+i-1)}).
\end{multline*}
\end{lem}
\begin{proof}
Consider the embeddings
\begin{equation}
\label{ment}
V_{(1^{2n+i})}\T \pi_j\hk V_{(1^{2n+i+1}(j+1))}\hk V_{(1^{2n+2+i})}\T \pi_j,
\end{equation}
where the first embedding comes from $(\ref{tens})$ and the second from
$(\ref{embV})$. We note that $(\ref{ment})$ means that
$\lim_{n\to\infty} V_{(1^{2n+i+1}(j+1))}\simeq L_{i,1}\T \pi_j$.
Note that
$$\widetilde \ch_q V^{2s+i}_{(1^{2n+i})}=q^{s(s+i)}\bin{2n+i}{n-s}.$$
Therefore,
$$\widetilde \ch_q V^{2s+i}_{(1^{2n+i})}-\ch_q L_{i,1}^{2s+i}=
q^{s(s+i)} \left(\bin{2n+i}{n-s}-\frac{1}{(q)_\infty}\right)=
O(q^{n+s(s+i-1)}).$$
We obtain that
$$\widetilde \ch_q (V_{(1^{2n+i+1})}\T\pi_j)^{2s+i+j}-
\ch_q (L_{i,1}\T\pi_j)^{2s+i+j}=O(q^{n+s(s+i-1)}).$$
To finish the proof it suffices to use $(\ref{ment})$ and the formula
$$\widetilde K_{2s+i+j, (1^{2n+i+1}(j+1))}(q)=
\widetilde \ch_q V_{(1^{2n+i+1}(j+1))}^{2s+i+j}
-\widetilde \ch_q V_{(1^{2n+i+1}(j+1))}^{2s+i+j+2}.$$
\end{proof}

We derive from this lemma that in $(\ref{chN})$ we can replace
reversed Kostka polynomials by the difference of the
characters of the weight subspaces of $L_{i,1}$.

Let $i=0$. Then for such $l$ that $j+l$ is even we obtain
\begin{multline*}
\lim_{n\to\infty} \ch N_{l,(1^{2n+i-1} (j+1))}=
q^{\frac{j(j+2)}{4(k+2)}-\frac{l(l+2)}{4(k+3)}}\times\\ \biggl(
\sum_{p\ge 0} q^{-(k+3)p^2-(l+1)p}
\frac{1}{(q)_\infty}(q^{(k+3)p+(l-j)/2}-q^{(k+3)p+(l+j+2)/2})-\\
\sum_{p> 0} q^{-(k+3)p^2+(l+1)p}
\frac{1}{(q)_\infty}(q^{(k+3)p-(l+j+2)/2}-q^{(k+3)p+(j-l-2)/2})\biggr)=
\end{multline*}
\begin{multline*}
\frac{q^{\frac{j(j+2)}{4(k+2)}-\frac{l(l+2)}{4(k+3)}+\frac{(l-j)^2}{4}}}
{(q)_\infty}
\biggl(\sum_{p\in\Z} q^{p^2(k+2)(k+3)+p((k+3)(j+1)-(k+2)(l+1))}-\\
\sum_{p\in\Z} q^{p^2(k+2)(k+3)+p((k+3)(j+1)+(k+2)(l+1))+(l+1)(j+1)}\biggr)=\\
q^{\triangle_{j+1,l+1}}
\biggl(\sum_{p\in\Z} q^{p^2(k+2)(k+3)+p((k+3)(j+1)-(k+2)(l+1))}-\\
\sum_{p\in\Z} q^{p^2(k+2)(k+3)+p((k+3)(j+1)+(k+2)(l+1))+(l+1)(j+1)}\biggr)=
\ch M_{j+1,l+1}(k+2,k+3)
\end{multline*}

One can repeat the same computation for $i=1$. We obtain
$$L_{i,1}\T L_{j,k}=\bigoplus_{l=0}^{k+1} L_{l,k+1}\T
M_{j+1,l+1}(k+2,k+3),$$
where the sum is taken over $l$ such that $l+i+j$ is even.

\subsection{The fermionic formula.}
Recall (see \cite{ss,kos}) that
\begin{equation}
\label{ff}
K^{(k)}_{l,\bm}(q)=
\sum_{\genfrac{}{}{0pt}{}{\bs\in\Z^k_{\ge 0}}{2|\bs|=|\bm|-l}}
q^{\bs A \bs+\bv \bs} \bin{A(\bm-2\bs)-\bv+\bs}{\bs},
\end{equation}
where $A_{\al,\be}=\min(\al,\be)$, $\bv_\al=\max(0, \al-k+l)$,
$|\bm|=\sum_{\al=1}^k m_\al$ and for two vectors
$\bm, {\bf n}\in\Z_{\ge 0}^k$ we set
$\bin{\bm}{{\bf n}}=\prod_{\al=1}^k \bin{m_\al}{n_\al}$.
We now explain how this formula naturally appears as a $q$-multiplicity.

\begin{lem}
\label{coinv}
Denote by $e(z)^i_s$ the coefficient in front of the power $z^s$ in the
series $(\sum_{j\in\Z} e_j z^j)^i$. Then
$$\ch_q L_{\bm,k}/
\bra  e_s, s\le 0;\ e(z)^{k-l+1}_s, s\le k-l+1;\ h_0+l \ket=
\widetilde K^{(k)}_{l,\bm}(q).$$
\end{lem}

\begin{proof}
We recall  an embedding
$V_\bm\hk L_{\bm,k}$. Let $v_\bm$ be a lowest weight vector of
$V_\bm$. We consider the principal
subspace
$$W=\C[e_N,e_{N-1},\ld]\cdot v_\bm,$$
where $N$ is fixed by $e_N v_\bm\ne 0$ and $e_{N+1} v_\bm=0$.
It is proved in \cite{two} that  the defining relations in $W$ are
$e(z)^{k+1}v_\bm=0$
and
$$e(z)^i v_\bm= z^{Ni-im_1-\cd -m_i} q(z^{-1}),\ i=1,\cd ,k,$$
where $q$ is some series. This means that the dual space of the quotient
$$W/\bra e_s, s\le 0;\ e(z)^{k-l+1}_s, s\le k-l+1;\ h_0+l \ket$$ can be
identified with a subspace of symmetric polynomials satisfying the conditions
\begin{enumerate}
\item The number of variables is $s=\frac{1}{2} (\sum_{i=1}^k im_i -l)$,
\item $f(\underbrace{z,\ld ,z}_{k+1},z_{i+1},\ld, z_s)=0$,
\item  $\deg_z f(\underbrace{z,\ld ,z}_a,z_{i+1},\ld, z_s)\le
\sum_{i=1}^k \min(a,i) m_i -a$,
\item $f(0,z_2,\ld,z_s)=0$,
\item $f(\underbrace{z,\ld ,z}_{k-l+1},z_{i+1},\ld, z_s)\div z^{k-l+2}$.
\end{enumerate}
But this space of symmetric polynomials coincides with the dual space
$(V_\bm/\bra h_0+l, e_0, e_1^{k-l+1} \ket)^*$ from \cite{kos},
and  the character of the latter
coincides with the corresponding restricted Kostka polynomial.

We recall that
$L_{\bm,k}=\lim_{N\to\infty} V_{(m_1,\ld, m_{k-1}, m_k+2N)}$.
Now our lemma follows from the equality $(\ref{indep})$.
\end{proof}

\begin{lem}
\label{qmult}
The $q$-multiplicity of $L_{l,k}$ in the decomposition of $L_{\bm,k}$ is
equal to the character of
\begin{equation}
\label{bra}
L_{\bm,k}/
\bra  e_s, s\le 0;\ e(z)^{k-l+1}_s, s\le k-l+1;\ h_0+l \ket.
\end{equation}
\end{lem}

\begin{proof}
We first note that the character of the space of the highest weight vectors
in $L_{\bm,k}$ of the weight $l$
is less or equal then the character of the quotient $(\ref{bra})$. But their
dimensions coincide.
\end{proof}

It is shown in \cite{kos} that the character of the space of symmetric
polynomials with the conditions from Lemma $\ref{coinv}$ is given by
the fermionic formula $(\ref{ff})$. This means that the $q$-multiplicities
from Lemma $\ref{qmult}$ are given by the fermionic formula.

\subsection{The limit of the fermionic formula.}
We want to find  the limit
$$\lim_{N\to\infty} \widetilde K^{(k+1)}_{l,(1^N (j+1))}(q),\
N+1=l+j \mod 2.$$ Up to a power of $q$ this limit coincides with the
unitary Virasoro character.

\begin{lem}
\label{h(m)}
Let $p(\bm)=\#\{\al=1,\ld,k:\ m_\al+\cd +m_k \text{ is odd }\}.$
Define
$$h(\bm)=\max\{\deg_q v:\ v\in V_\bm\}.$$
Then
$h(\bm)=\frac{\bm A \bm-p(\bm)}{4}.$
\end{lem}
\begin{proof}
We use the embedding
$(\ref{embV})$. Note that $h(1^N)=\frac{N^2}{4}-\frac{p(1^N)}{4}$. Therefore,
it follows from $(\ref{embV})$ that
\begin{equation}
\label{ineq}
h(\bm)\le \frac{(m_1+\cd +m_k)^2}{4}+\cd +\frac{m_k^2}{4}-\frac{p(\bm)}{4}=
\frac{\bm A \bm -p(\bm)}{4}.
\end{equation}
To prove that in $(\ref{ineq})$ we have an equality it suffices to show that
$$v(i_1)\T\cd \T v(i_k)\in U(\slt\T\C[t])\cdot \left(v(m_1+\cd +m_k)\T
\cd\T v(m_k)\right)$$
($i_\al=0$ if $m_\al+\cd +m_k$ is even and $i_\al=1$ otherwise).
But this follows from the formula
\begin{multline*}
\frac{1}{s!} e_{m_1+\cd +m_k-1}^s (v(m_1+\cd +m_k)\T \cd\T v(m_k))=\\
v(m_1+\cd +m_k-2)\T \cd\T v(m_s+\cd +m_k-2)\T v(m_{s+1}+\cd +m_k)\T\cd \T
v(m_k),
\end{multline*}
where $s$ is determined by $m_1=\cd =m_{s-1}=0$ and $m_s\ne 0$.
\end{proof}

In view of Lemma $\ref{h(m)}$ and a fact
${\genfrac{[}{]}{0pt}{0}{a}{b}}_{q^{-1}}=q^{-b(a-b)}
{\genfrac{[}{]}{0pt}{0}{a}{b}}_q$
we get:
\begin{multline*}
q^{\frac{p(\bm)}{4}}\widetilde K^{(k+1)}_{l,\bm}(q)=
q^{\frac{\bm A\bm}{4}}
\sum_{\genfrac{}{}{0pt}{}{\bs\in\Z^k_{\ge 0}}{2|\bs|=|\bm|-l}}
q^{-\bs A \bs-\bv \bs}
{\genfrac{[}{]}{0pt}{0}{A(\bm-2\bs)-\bv+\bs}{\bs}}_{q^{-1}}=\\
\sum_{\genfrac{}{}{0pt}{}{\bs\in\Z^k_{\ge 0}}{2|\bs|=|\bm|-l}}
q^{\frac{\bm A\bm}{4}-\bs A \bs-\bv \bs-\bs(A(\bm-2\bs)-\bv+\bs))}
\bin{A(\bm-\bs)-\bv+\bs}{\bs}=\\
\sum_{\genfrac{}{}{0pt}{}{\bs\in\Z^k_{\ge 0}}{2|\bs|=|\bm|-l}}
q^{(\frac{\bm}{2}-\bs) A (\frac{\bm}{2}-\bs)}
\bin{A(\bm-2\bs)-\bv+\bs}{\bs}.
\end{multline*}
Now let $\bm=(1^N(j+1))\in\Z^{k+1}$, $N+1=j+l\mod 2$.
We first rewrite the power of $q$ in the last line of the above formula
using the relation $2|\bs|=N+j+1-l$:
%\begin{multline*}
$$
(\frac{\bm}{2}-\bs) A (\frac{\bm}{2}-\bs)=
\sum_{\al,\be=1}^{k+1} s_\al s_\be +\frac{1}{4} (N^2+j+1+2N)-
N \sum_{\al=1}^{k+1} s_\al-  \sum_{\al=1}^{k+1} \min(\al,j+1)s_\al.
$$
%\end{multline*}
Replacing $s_1$ by $\frac{N+j+1-l}{2}-\sum_{\al=2}^{k+1} \al s_\al$ we get
\begin{multline}
\label{power}
\sum_{\al,\be=2}^{k+1} \max(\al,\be)(\min(\al,\be)-1) s_\al s_\be + \\
\sum_{\al=2}^{k+1} s_\al ((j+1-l)(\al-1)+\max(0,\al-j-1))+
\frac{(l-j)^2+j}{4}.
\end{multline}

Now we consider the binomial coefficient  $\bin{A(\bm-2\bs)-\bv+\bs}{\bs}$.
This is the product
$$\prod_{\al=1}^{k+1}
\bin{N+\min(\al,j+1)-2\sum_{\be=1}^{k+1} \min(\al,\be) s_\be -v_\al+s_\al}
{s_\al}.$$
Let $\al>1$. Then
\begin{multline}
\label{bin}
\bin{N+\min(\al,j+1)-2s_1-2\sum_{\be=2}^{k+1} \min(\al,\be) s_\be -v_\al+
s_\al} {s_\al}=\\
\bin{2\sum_{\be=2}^{k+1} (\be-\min(\al,\be)) s_\be +l-j-1+\min(\al,j+1)
-v_\al+ s_\al} {s_\al}.
\end{multline}
Now let $\al=1$. In this case the binomial coefficient
depends  on $N$. We want to know the limit of this expression while
$N\to\infty$.
\begin{multline}
\label{fac}
\bin{N+1-s_1-2\sum_{\be=2}^{k+1} s_\be -v_1} {s_1}=\\
\bin{\frac{N+l+1-j}{2}+ \sum_{\be=2}^{k+1} (\be-2) s_\be -v_1}
{2\sum_{\be\ge 2} (\be-1)s_\be+l-j -v_1}\to
\frac{1}{\fac{2\sum_{\be\ge 2} (\be-1)s_\be+l-j -v_1}}.
\end{multline}

We obtain the following proposition.
\begin{prop}
Fix $0\le j\le k$, $0\le l\le k+1$. Then
\begin{multline*}
\ch M_{j+1,l+1}(k+2,k+3)=
q^{\triangle_{j+1,l+1}}
\sum_{\bt=(t_2,\ld,t_{k+1})\in\Z^k_{\ge 0}} q^{\bt B \bt +\bu\bt}\\
\frac{\prod\limits_{\al=2}^{k+1}
\bin{2\sum_{\be=2}^{k+1} (\be-\min(\al,\be)) t_\be +l-j-1+\min(\al,j+1)
-v_\al+ t_\al} {t_\al}
}
{\fac{2\sum_{\be\ge 2} (\be-1)t_\be+l-j -v_1}},
\end{multline*}
where
$B_{\al,\be}=\max(\al,\be)(\min(\al,\be)-1)$,
$\bu_\al=(j+1-l)(\al-1)+\max(0,\al-j-1)$.
\end{prop}
\begin{proof}
In view of a formulas $(\ref{power}),(\ref{bin})$ and $(\ref{fac})$ we
only need to find the power of $q$ in front of the sum. Note that
$p(1^N (j+1))=j+i$, where $i=0,1$ and $i+N+1$ is even. Now it is
enough to mention that
\begin{multline*}
\frac{i(i+1)}{12}+\frac{j(j+2)}{4(k+2)}-\frac{l(l+2)}{4(k+3)}-
\frac{i+j}{4}+\frac{(l-j)^2+j}{4}=\\
\frac{((k+3)j-(k+2)l+1)^2-1}{4(k+2)(k+3)}=\triangle_{j+1,l+1}(k+2,k+3).
\end{multline*}
(The first three terms come from the action of the Sugawara operators, the
fourth one from $p(\bm)$ and the last from the formula $(\ref{power})$).
\end{proof}

\appendix
\section{The ABF finitization and Kostka polynomials}
In this appendix we study the connection between ABF finitization
of the minimal Virasoro unitary characters (see \cite{ABF}) and Kostka
polynomials $K^{(k)}_{l,(1^N (j+1))}$, generalizing the case $j=0$ in
\cite{JMT}.

We first recall the ABF finitization. Fix some $a,b\in\Z$. For $N\ge 0$ such
that $N\equiv b-a\mod 2$, define a polynomial
\begin{multline*}
\widehat\chi_{b,a}^{(r,r+1)}(q;N)=\sum_{n\in\Z}
q^{r(r+1)n^2+((r+1)b-ra)n}\bin{N}{\frac{N-b+a}{2}-(r+1)n}\\
-\sum_{n\in\Z} q^{r(r+1)n^2+((r+1)b+ra)n+ba}\bin{N}{\frac{N-b-a}{2}-(r+1)n}.
\end{multline*}
In view of the Rocha-Caridi formula it is obvious that
$$\lim_{\genfrac{}{}{0pt}{}{N\to\infty}{N\equiv b-a \mod 2}}
q^{\triangle_{b,a}} \widehat\chi_{b,a}^{(r,r+1)}(q;N)=\ch M_{b,a}(r,r+1).$$

\begin{lem}
\begin{multline*}
q^{\frac{N^2}{4}-\frac{(a-b)^2}{4}}\widehat\chi_{b,a}^{(r,r+1)}(q^{-1};N)=
\sum_{n\in\Z}
q^{(r+1)n^2-na}\bin{N}{\frac{N-b+a}{2}-(r+1)n}\\
-\sum_{n\in\Z}q^{(r+1)n^2+na}\bin{N}{\frac{N-b-a}{2}-(r+1)n}.
\end{multline*}
\end{lem}
\begin{proof}
We just use the identity
$$\genfrac{[}{]}{0pt}{}{N}{k}_{q^{-1}}=q^{-k(N-k)}\bin{N}{k}.$$
\end{proof}

\begin{prop}
\begin{multline*}
q^{\frac{(l-j)^2}{4}} K^{(k)}_{l,(1^N(j+1))}(q)=\\
q^{\frac{(N+1)^2}{4}}\sum_{s=0}^j
\widehat\chi_{j+1-2s,l+1}^{(k+1,k+2)}(q^{-1};N+1)-
q^{\frac{N^2}{4}}\sum_{s=0}^{j-1}
\widehat\chi_{j-2s,l+1}^{(k+1,k+2)}(q^{-1};N).
\end{multline*}
\end{prop}

\begin{proof}
In view of Corollary $\ref{fk}$ and the alternating sum formula we have
\begin{multline*}
K^{(k)}_{l,\bm}(q)=\\
\sum_{p\ge 0} q^{(k+2)p^2+(l+1)p} K_{2(k+2)p+l,\bm}(q)-
\sum_{p> 0} q^{(k+2)p^2-(l+1)p} K_{2(k+2)p-l-2,\bm}(q)=\\
\sum_{p\in\Z} q^{(k+2)p^2+(l+1)p} (\ch_q V_\bm^{2(k+2)p+l}-
\ch_q V_\bm^{2(k+2)p+l+2}).
\end{multline*}
Now let $\bm=(1^N(j+1))$. We show by induction on $j$ that
\begin{equation}
\label{chj}
\ch_q V_{(1^N(j+1))}^l=\bin{N+1}{\frac{N-j+1-l}{2}}+
\sum_{s=0}^{j-1}
q^{\frac{N+j+1-l-2s}{2}}\bin{N}{\frac{N+j+1-l-2s}{2}}.
\end{equation}
For $j=0$ we have $\ch_q V_{(1^{N+1})}^l=\bin{N+1}{\frac{N+1-l}{2}}.$
For the induction procedure we use $(\ref{dem})$. This gives
\begin{multline*}
\ch_q V_{(1^N(j+1))}^l=\ch_q V_{(1^N j)}^{l+1}+
q^{\frac{N+j+1-l}{2}} \ch_q V_{(1^N)}^{l-j-1}=\\
\bin{N+1}{\frac{N-j+1-l}{2}}+
\sum_{s=0}^{j-1}
q^{\frac{N+j+1-l-2s}{2}}\bin{N}{\frac{N+j+1-l-2s}{2}}.
\end{multline*}
Using the identity
$$q^b\bin{a}{b}+\bin{a}{b-1}=\bin{a+1}{b}$$ we rewrite  $(\ref{chj})$
in the following way
$$\ch_q V_{(1^N(j+1))}^l=\sum_{s=0}^j \bin{N+1}{\frac{N+j+1-l-2s}{2}}-
\sum_{s=0}^{j-1} \bin{N}{\frac{N+j-1-l-2s}{2}}.$$
Therefore, we obtain
\begin{multline*}
K^{(k)}_{l,(1^N(j+1))}(q)=
\sum_{p\in\Z} q^{(k+2)p^2+(l+1)p}\times \\
\sum_{s=0}^j \left( \bin{N+1}{\frac{N+j+1-l-2s}{2}-(k+2)p}-
             \bin{N+1}{\frac{N+j-1-l-2s}{2}-(k+2)p}\right) -\\
\shoveleft{\sum_{p\in\Z} q^{(k+2)p^2+(l+1)p}}\times\\
\sum_{s=0}^{j-1} \left(\bin{N}{\frac{N+j-1-l-2s}{2}-(k+2)p}-
             \bin{N}{\frac{N+j-3-l-2s}{2}-(k+2)p}\right).
\end{multline*}
We show that
\begin{multline}
\label{N+1}
\sum_{p\in\Z} q^{(k+2)p^2+(l+1)p}
\Biggl(\bin{N+1}{\frac{N+j+1-l-2s}{2}-(k+2)p}-\\
             \bin{N+1}{\frac{N+j-1-l-2(j-s)}{2}-(k+2)p}\Biggr)=
q^{\frac{(N+1)^2-(l-j)^2}{4}}
\widehat\chi_{j+1-2s,l+1}^{(k+1,k+2)}(q^{-1};N+1).
\end{multline}
We rewrite the left-hand side of $(\ref{N+1})$ as
\begin{multline*}
\sum_{p\in\Z} q^{(k+2)p^2-(l+1)p}
\bin{N+1}{\frac{N+j+1-l-2s}{2}+(k+2)p}-\\
\sum_{p\in\Z} q^{(k+2)p^2+(l+1)p}
\bin{N+1}{\frac{N+1-(l+1)-(j+1-2s)}{2}-(k+2)p}=
\end{multline*}
\begin{multline*}
=\sum_{p\in\Z} q^{(k+2)p^2-(l+1)p}
\bin{N+1}{\frac{N+1+(l+1)-(j+1-2s)}{2}-(k+2)p}-\\
\sum_{p\in\Z} q^{(k+2)p^2+(l+1)p}
\bin{N+1}{\frac{N+1-(l+1)-(j+1-2s)}{2}-(k+2)p},
\end{multline*}
which coincides with the right-hand side of $(\ref{N+1})$.
In the same way one can prove that
\begin{multline*}
\sum_{p\in\Z} q^{(k+2)p^2+(l+1)p}
\Biggl(\bin{N}{\frac{N+j-1-l-2s}{2}-(k+2)p}-\\
             \bin{N}{\frac{N+j-3-l-2(j-1-s)}{2}-(k+2)p}\Biggr)=
q^{\frac{N^2-(l-j)^2}{4}}
\widehat\chi_{j-2s,l+1}^{(k+1,k+2)}(q^{-1};N).
\end{multline*}
Proposition is proved.
\end{proof}

\newcounter{a}
\setcounter{a}{2}

\end{document}